\documentclass{amsart}
\usepackage{amssymb}
\begin{document}
\newtheorem{theorem}{Theorem}[section]
\newtheorem{lemma}[theorem]{Lemma}
\def\operatorname#1{{\rm#1\,}}
\def\range{\operatorname{range}}
\def\Pspan{\operatorname{span}}
\def\RR{\mathcal{R}}
\def\rank{\operatorname{rank}}
\def\id{\operatorname{Id}}
\def\trace{\operatorname{trace}}
\def\imag{\operatorname{Im}}
\font\pbglie=eufm10
\def\SS{{\text{\pbglie S}}}
\def\spec{\operatorname{Spec}}
 \makeatletter
  \renewcommand{\theequation}{%
   \thesection.\arabic{equation}}
  \@addtoreset{equation}{section}
 \makeatother
\title{Szab\'o Osserman IP Pseudo-Riemannian manifolds}
\author{Peter B. Gilkey, Raina
Ivanova, and Tan Zhang}
\begin{address}{PG: Mathematics Department, University of Oregon, Eugene OR 97403 USA}
\end{address}
\begin{email}
{gilkey@darkwing.uoregon.edu}\end{email}
\begin{address}{RI: Mathmatics Department,
University of Hawaii - Hilo,
200 W. Kawili St.,
Hilo,  HI 96720 USA}\end{address}
\begin{email}
{rivanova@hawaii.edu}
\end{email}
\begin{address}{TZ: Department of Mathematics and Statistics,
Murray State University, Murray, KY 42071 USA}\end{address}
\begin{email}{tan.zhang@murraystate.edu}
\end{email}
\date{25 April 2002 Version 3j}
\begin{abstract} We construct a family of pseudo-Riemannian manifolds so that the
skew-symmetric curvature operator, the Jacobi operator, and the  Szab\'o operator have constant eigenvalues on their
domains of definition. This provides new and non-trivial examples of Osserman, Szab\'o, and IP manifolds. We also study when the
associated Jordan normal form of these operators is constant.
  {\it
Subject Classification}: 53B20.
\end{abstract}
\maketitle
\font\pbglie=eufm10
\def\Gr{\text{Gr}}
\section{Introduction}
Let $(M,g)$ be a pseudo-Riemannian manifold of signature $(p,q)$. Let $R$ be the Riemann curvature:
\begin{eqnarray*}
&&R(Z_1,Z_2):=\nabla_{Z_1}\nabla_{Z_2}-\nabla_{Z_2}\nabla_{Z_1}-\nabla_{[Z_1,Z_2]},\\
&&R(Z_1,Z_2,Z_3,Z_4):=(R(Z_1,Z_2)Z_3,Z_4).\end{eqnarray*}
We can use $R$ and $\nabla R$ to define several natural operators:\begin{enumerate}
\item The {\it Jacobi operator} $J(X):Y\rightarrow R(Y,X)X$ is a symmetric operator with $J(X)X=0$. It plays an
important role in the study of geodesic sprays. Since $J(cX)=c^2J(X)$, the natural domains of
definition for
$J$ are the pseudo-sphere bundles
$S^\pm(M,g):=\{X\in TM:(X,X)=\pm1\}$.
\item The {\it Szab\'o operator} $\SS(X):Y\rightarrow\nabla_XR(Y,X)X$ is a symmetric operator with $\SS(X)X=0$. It
plays an important role in the study of totally isotropic manifolds. Since $\SS(cX)=c^3\SS(X)$, the natural domains of
definition for $\SS$ are $S^\pm(M,g)$.
\item Let $\{X_1,X_2\}$ be an oriented orthonormal basis for a non-degenerate $2$ plane $\pi$. The {\it skew-symmetric curvature
operator}
$\RR(\pi):Y\rightarrow R(X_1,X_2)Y$  depends on the orientation of $\pi$ but not on the particular orthonormal
basis chosen. The natural domains of definition for $\RR(\cdot)$ are the oriented Grassmannians of timelike, mixed
(signature $(1,1)$), and spacelike $2$ planes.
\end{enumerate}

The {\it spectrum} $\spec(T)\subset\mathbb{C}$ of a linear map $T$ is the set of complex eigenvalues of $T$. It is natural
to ask what are the geometric constraints that are imposed by assuming that the spectrum (or more generally the Jordan normal form)
of one of these 3 natural operators is constant on the appropriate domains of definition. 

Let $(M,g)$ be a pseudo-Riemannian manifold.
$(M,g)$ is said to be {\it spacelike Osserman} if $\spec(J(\cdot))$ is constant
on $S^+(M,g)$,  $(M,g)$ is said to be {\it spacelike Szab\'o} if $\spec(\SS(\cdot))$ is constant
on $S^+(M,g)$, and  $(M,g)$ is said to be {\it spacelike IP}  if $\spec(\RR(\pi))$ is
constant on the Grassmannian of oriented spacelike $2$ planes in the tangent bundle $TM$. One defines {\it timelike Osserman},
{\it timelike Szab\'o}, {\it timelike IP}, and {\it mixed IP} similarly.
The eigenvalue $\{0\}$ plays a distinguished role. We say $(M,g)$ is {\it nilpotent Osserman} if $\spec(J(X))=\{0\}$
for all $X$, {\it nilpotent Szab\'o} and {\it nilpotent IP} are defined similarly.

 The names Osserman, Szab\'o, and IP are used because the seminal papers for
this subject in the Riemannian setting are due to Osserman \cite{refOss} for the operator $J(\cdot)$, to Szab\'o
\cite{refSzabo} for the operator $\SS(\cdot)$, to Ivanov and Petrova \cite{refIP} and Stanilov and Ivanova \cite{{refRIGS}} for
the operator
$\RR(\cdot)$. The spectral properties of the operators $J(\cdot)$ and $\RR(\cdot)$ have been studied extensively; we refer to
\cite{refGRKVL,refGil} for a more complete historical discussion and bibliography. By contrast, the operator
$\SS(\cdot)$ has received considerably less attention.

 Suppose $p\ge1$ and $q\ge1$. Then the notions
spacelike Osserman (resp. spacelike Szab\'o) and timelike Osserman ({resp. timelike Szab\'o}) are equivalent,
so one simply says that
$(M,g)$ is {\it Osserman} (resp. {\it Szab\'o}). Similarly, if $p\ge2$ and if $q\ge2$, then spacelike, mixed, and
timelike IP are equivalent notions so $(M,g)$ is said to be {\it IP}. See \cite{refGil} for details. 
We shall use the
words `nilpotent', `Osserman', `Szab\'o', and `IP' as adjectives. Thus, for example, to say that
a manifold is nilpotent Osserman Szab\'o IP means that it is simultaneously nilpotent Szab\'o, nilpotent Osserman, and
nilpotent IP.
We say that $(M,g)$ is {\it locally symmetric} if $\nabla R=0$ and {\it locally homogeneous} if the local isometries
of $(M,g)$ act transitively on $M$. We say that $(M,g)$ is {\it Ricci flat} if the Ricci tensor vanishes identically.

Let $( x, y)=(x_1,...,x_p,y_1,...,y_p)$ be coordinates on the manifold
$M:=\mathbb{R}^{2p}$. Let $\psi_{ij}(x)=\psi_{ji}(x)$ be a symmetric $2$ tensor $\psi$ on $\mathbb{R}^p$. We
define a non-degenerate pseudo-Riemannian metric $g_\psi$ of balanced signature $(p,p)$ on $M$ by setting:
\begin{equation}
ds^2_\psi:=\textstyle\sum_i dx^i\circ dy^i+\textstyle\sum_{i,j}\psi_{ij}(x)dx^i\circ dx^j.
\label{Eqn1.1}\end{equation}

In Section \ref{Sect2}, we will prove the following result:

\begin{theorem}\label{Thm1.1} Let $p\ge2$ and let $\psi(x)$ be a symmetric $2$ tensor. Then $(M,g_\psi)$ is:\begin{enumerate}
\item a pseudo-Riemannian manifold of signature $(p,p)$;\item nilpotent Szab\'o Osserman IP;
\item Ricci flat and Einstein;
\item neither locally homogeneous nor locally symmetric for generic $\psi$.\end{enumerate}
\end{theorem}

Nilpotent Osserman manifolds have been constructed previously \cite{refBBGZ,refGRKVL}. We can describe one family which arises from
affine geometry as follows. Let
$\Gamma_{ij}{}^k(x)$ be the Christoffel symbols of an arbitrary torsion free connection
$\nabla$ on $\mathbb{R}^p$. Let $R_\nabla$ be the associated curvature operator and let 
$J_\nabla(X_1):X_2\rightarrow R_\nabla(X_2,X_1)X_1$ be
the associated Jacobi operator on $\mathbb{R}^p$. We say that the connection $\nabla$ is {\it nilpotent affine Osserman} if for all
$X\in T\mathbb{R}^p$, $\spec(J_\nabla(X))=\{0\}$,  i.e. $J_\nabla(X)^p=0$.

Following  Garci\'a-Ri\'o, Kupeli, and V\'azquez-Lorenzo
\cite{refGRKVL} (see page 147), define an associated metric on $M=\mathbb{R}^{2p}$ by setting:
\begin{equation}ds^2_\nabla=\textstyle\sum_idx^i\circ dy^i-2\sum_{ijk}y_k\Gamma_{ij}{}^k(x)dx^i\circ dx^j.
\label{Eqn112}\end{equation}
 Then
$(M,ds^2_\nabla)$ is nilpotent Osserman if and only if
$\nabla$ is nilpotent affine Osserman. This metric is quite different in flavor from ours
as the coefficients depend on the $y$ variables as well as on the $x$ variables. There does not
seem to be any direct connection between the metrics defined in equations (\ref{Eqn1.1}) and (\ref{Eqn112}).

Nilpotent IP manifolds have also been constructed previously \cite{refGiZa}.
However, comparatively little is known about Szab\'o manifolds - see
\cite{refGIZ,refGilStav} for some preliminary results in the algebraic setting. In particular, the manifolds $(M,g_\psi)$ are the
only known irreducible Szab\'o manifolds which are not locally symmetric.

\medbreak The eigenvalue structure does not determine the conjugacy class (i.e. the real Jordan normal form) of a
symmetric or skew-symmetric linear operator in the higher signature setting. We
will use the words `timelike', `spacelike', and `Jordan' as adjectives. Thus, for example, to say $(M,g)$ is timelike
Jordan IP means that the Jordan normal form of the skew-symmetric curvature operator is constant on the
Grassmannian of timelike oriented $2$-planes.  We shall omit the accompanying adjectives `timelike and spacelike' if
both apply. Thus
$(M,g)$ is Jordan Osserman means
$(M,g)$ is both timelike Jordan Osserman and spacelike Jordan Osserman, i.e. the Jordan normal form of $J$ is
constant on $S^+(M,g)$ and on $S^-(M,g)$.

There are no known timelike or spacelike Jordan Szab\'o manifolds which are not locally symmetric. Section \ref{Sect3} is devoted to
the proof of:
\begin{theorem}\label{Thm1.2} If $(M,g_\psi)$ is not locally symmetric, then $(M,g_\psi)$ is neither spacelike Jordan Szab\'o nor
timelike Jordan Szab\'o.\end{theorem}

It is useful to consider a subfamily of the metrics defined in equation (\ref{Eqn1.1}). Let $f$ be a real-valued function
on
$\mathbb{R}^p$. In equation (\ref{Eqn1.1}) we set
$\psi=df\circ df$ and define
\begin{equation}ds_f^2:=\textstyle\sum_idx^i\circ dy^i+\textstyle\sum_{i,j}\frac{\partial f}{\partial x_i}\frac{\partial f}{\partial
x_j} dx^i\circ dx^j.\label{Eqn1.2}\end{equation}
We can realize $(M,g_f)$ as a hypersurface in a flat space. Let
$\{\alpha_1,...,\alpha_p,\beta_1,...,\beta_p,\gamma\}$ and
$\{\alpha_1^*,...,\alpha_p^*,\beta_1^*,...,\beta_p^*,\gamma^*\}$
be a basis and associated dual basis for vector spaces $W$ and $W^*$, respectively. We give $W$
the inner product of signature $(p,p+1)$:
$$ds^2_W:=\alpha_1^*\circ\beta_1^*+...+\alpha_p^*\circ\beta_p^*+\gamma^*\circ\gamma^*.$$
Let $F:M\rightarrow W$ be the isometric embedding:
$$F(x, y):=\textstyle\sum_i(x_i\alpha_i+y_i\beta_i)+f( x)\gamma.$$
The normal $\nu$ to the hypersurface is given by
$\nu=-\textstyle\frac{\partial f}{\partial x_i}\beta_i+\gamma$,
so the second fundamental form of the embedding is
\begin{equation}L(Z_1,Z_2)=Z_1Z_2(f).\label{Eqn1.3}\end{equation}
We define distributions
$$\mathcal{X}:=\Pspan\{\partial_1^x,...,\partial_p^x\}\quad\text{ and }\quad
  \mathcal{Y}:=\Pspan\{\partial_1^y,...,\partial_p^y\}.$$
We then have $L(Z_1,Z_2)=0$ if $Z_1\in\mathcal{Y}$ or $Z_2\in\mathcal{Y}$ so the restriction $L_{\mathcal{X}}$ of the second
fundamental form
$L$ to the distribution $\mathcal{X}$ carries the essential information. If, for
example, we set $f(x)=\sum_i\varepsilon_ix_i^2$, then:
$$L_{\mathcal{X}}(\partial_i^x,\partial_j^x)=\varepsilon_i\delta_{ij}$$
so there are non-trivial examples where $L_{\mathcal{X}}$ is non-degenerate on
$M$. In Sections \ref{Sect4} and \ref{Sect5}, we prove:

\begin{theorem}\label{Thm1.3} Assume the quadratic form $L_{\mathcal X}$ is non-degenerate. \begin{enumerate}\item If $p=2$ or if 
$p\ge3$ and if
$L_{\mathcal{X}}$ is definite, then $(M,g_f)$ is:
\begin{enumerate}\item nilpotent
Jordan Osserman;
\item nilpotent spacelike and timelike Jordan IP;
\item not mixed Jordan IP.\end{enumerate}
\item If $p\ge3$ and if $L_{\mathcal X}$ is indefinite, then $(M,g_f)$ is:
\begin{enumerate}\item neither
 spacelike Jordan Osserman nor timelike Jordan Osserman;
\item  nilpotent spacelike and timelike
Jordan IP;
\item not mixed IP.\end{enumerate}\end{enumerate}\end{theorem}

So far, we have discussed the balanced (or neutral) signature $p=q$. There are also results available when $p\ne q$. Give
$\mathbb{R}^{(a,b)}$ the canonical flat metric of signature $(a,b)$. We conclude the paper in Section \ref{Sect6} by proving:

\goodbreak\begin{theorem}\label{Thm1.5} Assume the quadratic form $L_{\mathcal{X}}$ of $(M,g_f)$ is positive definite.
Let $N:=M\times\mathbb{R}^{(a,b)}$ and let $g_N$ be the product metric on $N$.
\begin{enumerate}
\item $(N,g_N)$ is a nilpotent Osserman Szab\'o IP manifold of signature $(p+a,q+b)$.
\item For generic $f$, $(N,g_N)$ is
\begin{enumerate}
\item neither spacelike Jordan Szab\'o,
\item nor timelike Jordan Szab\'o,
\item nor locally
homogeneous,
\item nor locally symmetric.\end{enumerate}
\item $(N,g_N)$ is not mixed Jordan IP.
\item Suppose that $b=0$. Then $(N,g_N)$ is:
\begin{enumerate}
\item neither timelike Jordan Osserman nor timelike Jordan IP;
\item spacelike Jordan Osserman and spacelike Jordan IP.
\end{enumerate}
\item Suppose that $a=0$. Then $(N,g_N)$ is:
\begin{enumerate}\item timelike Jordan Osserman and timelike Jordan IP;
\item neither spacelike Jordan Osserman nor spacelike Jordan IP.
\end{enumerate}
\item Suppose that $a>0$ and $b>0$. Then $(N,g_N)$ is:
\begin{enumerate}\item  neither timelike Jordan Osserman nor timelike Jordan IP;
\item neither spacelike Jordan
Osserman nor spacelike Jordan IP.
\end{enumerate}
\end{enumerate}\end{theorem}

\medbreak\noindent{\bf Note:} By Theorem \ref{Thm1.3},  Jordan Osserman and Jordan IP are different notions. By Theorem
\ref{Thm1.5}, timelike Jordan Osserman (resp. timelike Jordan IP) and spacelike Jordan Osserman (resp. spacelike Jordan
IP) are different notions as well.

\section{Nilpotent Jordan Szab\'o Osserman manifolds}\label{Sect2}

We begin the proof of Theorem \ref{Thm1.1} by determining the
curvature tensor of
$(M,g_\psi)$. Let
$\psi_{ij/k}:=\partial_k^x\psi_{ij}$ and let $\psi_{ij/kl}:=\partial_k^x\partial_l^x\psi_{ij}$.

\begin{lemma}\label{Lem2.1} Let $Z_\nu$ be vector fields on
$(M,g_\psi)$. We have:
\begin{enumerate}
\item $\nabla\partial_i^y=0$;
\item $R(Z_1,Z_2,Z_3,Z_4)=0$ if one of the $Z_\nu\in\mathcal{Y}$ for $1\le\nu\le4$;
\item $\nabla R(Z_1,Z_2,Z_3,Z_4;Z_5)=0$ if one of the $Z_\nu\in\mathcal{Y}$ for $1\le\nu\le5$;
\item $R(\partial_i^x,\partial_j^x,\partial_k^x,\partial_l^x)=-\frac12(\psi_{il/jk}+\psi_{jk/il}-\psi_{ik/jl}-\psi_{jl/ik})$.
\end{enumerate}\end{lemma}

\begin{proof} Let $Z_1$ and $Z_2$ be coordinate vector fields. We then have:
$$(\nabla_{Z_1}\partial_i^y,Z_2)=\textstyle{\frac12}\{\partial_i^y(Z_1,Z_2)+Z_1(\partial_i^y,Z_2)-Z_2(Z_1,\partial_i^y)\}=0.$$
Assertion (1) now follows. We use it to see
\begin{equation}R(Z_1,Z_2,\partial_i^y,Z_3)=
((\nabla_{Z_1}\nabla_{Z_2}-\nabla_{Z_2}\nabla_{Z_1}-\nabla_{[Z_1,Z_2]})\partial_i^y,Z_3)=0.
\label{Eqn2.1}\end{equation}
This proves assertion (2) if $Z_3\in\mathcal{Y}$; the curvature symmetries
then show that $R(Z_1,Z_2,Z_3,Z_4)=0$ if any of the remaining vectors belong to $\mathcal{Y}$. Since
$\nabla\partial_i^y=0$, we can
covariantly differentiate equation (\ref{Eqn2.1}) and get
\begin{equation}\nabla_{Z_1}R(Z_2,Z_3,\partial_i^y,Z_4)=0;
\label{eqn2.2RI}\end{equation}
assertion (3) now follows from equation (\ref{eqn2.2RI}) and from the curvature symmetries.

Since $(\nabla_{\partial_i^x}\partial_j^x,\partial_k^y)=0$ and 
$(\nabla_{\partial_i^x}\partial_j^x,\partial_k^x)=\textstyle\frac12(\psi_{ik/j}+\psi_{jk/i}-\psi_{ij/k})$, we have
$$\nabla_{\partial_i^x}\partial_j^x=\textstyle\frac12\sum_k(\psi_{ik/j}+\psi_{jk/i}-\psi_{ij/k})\partial_k^y.$$
We complete the proof of the Lemma by computing:
\begin{eqnarray*}
&&R(\partial_i^x,\partial_j^x,\partial_k^x,\partial_l^x)=((\nabla_{\partial_i^x}\nabla_{\partial_j^x}-\nabla_{\partial_j^x}\nabla_{\partial_i^x})\partial_k^x,\partial_l^x)\\
&=&\textstyle\frac12(\{\partial_i(\psi_{j\nu/k}+\psi_{k\nu/j}-\psi_{jk/\nu})-
  \partial_j(\psi_{i\nu/k}+\psi_{k\nu/i}-\psi_{ik/\nu})\}\partial_\nu^y,\partial_l^x)\\
&=&\textstyle\frac12(\psi_{jl/ki}+\psi_{kl/ji}-\psi_{jk/il}-\psi_{il/jk}-\psi_{kl/ij}+\psi_{ik/jl}).
\end{eqnarray*}
\end{proof}

\begin{proof}[Proof of Theorem \ref{Thm1.1}] 
\par 1) It is clear from equation (\ref{Eqn1.1}) that $(M,g_\psi)$ is a pseudo-Riemannian
manifold of signature
$(p,p)$. 
\par 2) We use Lemma
\ref{Lem2.1} to see that
\begin{equation}
\mathcal{Y}\subset\ker(J(Z_1)),\quad
\mathcal{Y}\subset\ker(R(Z_1,Z_2)),\quad
\mathcal{Y}\subset\ker(\SS(Z_1)).\label{Eqn2.3}\end{equation} 
Since $\mathcal{Y}$ is a
totally isotropic subspace of $T_PM$ of dimension $p$, $\mathcal{Y}^\perp=\mathcal{Y}$.  Since $(R(Z_1,Z_2)Z_3,\partial_i^y)=0$ and
$(\nabla_{Z_1}R(Z_2,Z_3)Z_4,\partial_i^y)=0$, we have:
\begin{equation}
\range(J(Z_1))\subset\mathcal{Y},\quad \range(R(Z_1,Z_2))\subset\mathcal{Y},\quad
\range(\SS(Z_1))\subset\mathcal{Y}.\label{Eqn2.4}\end{equation}
We use equations (\ref{Eqn2.3}) and (\ref{Eqn2.4}) to show
\begin{equation}J(Z_1)^2=0,\quad\SS(Z_1)^2=0,\quad\text{and}\quad R(Z_1,Z_2)^2=0.\label{Eqn2.5}\end{equation}
This shows that
$$\spec(J(Z_1))=\{0\},\quad\spec(\SS(Z_1))=\{0\},\quad\text{and}\quad\spec(R(Z_1,Z_2))=\{0\}$$
for any vector fields
$Z_\nu$. Consequently $(M,g_\psi)$ is nilpotent Jordan Szab\'o IP.

3) Let $\varrho$ be the Ricci tensor. Since $\varrho(Z,Z)=\trace(J(Z))$ and since $J(Z)^2=0$, $\varrho(Z,Z)=0$. We polarize to see
that
$\varrho\equiv0$. Thus $(M,g_\psi)$ is Ricci flat and Einstein.

4) Clearly, $(M,g_\psi)$ is
generically neither locally homogeneous nor locally symmetric.\end{proof}
\section{Jordan Szab\'o manifolds}\label{Sect3}

Theorem \ref{Thm1.2} will follow from the following Lemma:

\begin{lemma}\label{Lem3.1} Let $p\ge2$ and let $P\in M$. If $\nabla
R_\psi(P)$ does not vanish identically, then
$\rank(\SS(\cdot))$ is constant neither on $S^+(T_PM)$ nor on $S^-(T_PM)$. Thus $(M,g_f)$ is neither
spacelike Jordan Szab\'o nor timelike Jordan Szab\'o.
\end{lemma}

\begin{proof} 
Suppose $\rank(\SS(\cdot))=r>0$ is constant on $S^+(T_PM)$; the timelike case is
similar. Let
$\mathcal{V}^+$ be a maximal spacelike subspace of $T_PM$ and 
let $\mathcal{V}^-:=(\mathcal{V}^+)^\perp$ be the complementary timelike
subspace. Let $\rho^\pm$ be orthogonal projection on $\mathcal{V}^\pm$.  
If $Z\in S^+(\mathcal{V}^+)$, then we define: 
$$\check\SS(Z):=\rho^+\SS(Z)\rho^+.$$

We wish to show that $\rank\check
\SS(Z)=r$. Let $\{Z_1,...,Z_r\}$ be tangent vectors at $P$ so
$\{\SS(Z)Z_1,...,\SS(Z)Z_r\}$ is a basis for $\range(\SS(Z))$. As
$T_PM=\mathcal{Y}+\mathcal{V}^+$, we may decompose $Z_i=V_i^++Y_i$, 
where $V_i^+\in \mathcal{V}^+$ and $Y_i\in \mathcal{Y}$. Since
$\mathcal{Y}\subset\ker\SS(Z)$,
$\SS(Z)Z_i=\SS(Z)V_i^+$ and thus $\{\SS(Z)V_1^+,...,\SS(Z)V_r^+\}$ is a basis
for $\range(\SS(Z))$. As
$\ker\rho^+=\mathcal{V}^-$ is timelike, as
$\mathcal{Y}$ is totally isotropic, and as $\range(\SS(Z))\subset \mathcal{Y}$, the vectors
$$\{\rho^+\SS(Z)\rho^+V_1^+,...,\rho^+\SS(Z)\rho^+V_r^+\}$$
are linearly independent. Consequently, $\rank(\check\SS(Z))\ge
r$. Since the reverse inequality is immediate, we have as desired that
$$\rank(\check\SS(Z))=r\quad\text{for}\quad Z\in S^{p-1}:=S^+(\mathcal{V}^+).$$

Since $\SS$ is self-adjoint and $\rho^+$ is self-adjoint,
$\check\SS(Z)$ is a self-adjoint map of $\mathcal{V}^+$. Let
$E_+$ and
$E_-$ be the span of the eigenvectors with positive and negative eigenvalues respectively; these are non-trivial as $r>0$. Since
$\check\SS$ is a self-adjoint map with constant rank, $E_+$ and $E_-$ are vector bundles over $S^{p-1}$. Since $\check\SS$
is self-adjoint, and since $\check\SS(Z)Z=0$,
$$E_+(Z)\perp E_-(Z)\quad\text{and}\quad E_\pm(Z)\subset Z^\perp=T_ZS^{p-1}.$$

Let $N$ be the north pole of $S^{p-1}$. Since $S^{p-1}-\{N\}$ is contractable, there exists a section $s_+$ to $E_+$
vanishing only at $N$. Since $\SS(-Z)=-\SS(Z)$, $E_+(Z)=E_-(-Z)$. Thus
$s_-(Z):=s_+(-Z)$ is a section to $E_-$ which only vanishes at $-N$. Since
$E_+\perp E_-$, we have $s_+(Z)\perp s_-(Z)$. Consequently, the vector field
$$s(Z):=s_+(Z)+s_-(Z)$$
is nowhere vanishing on $S^{p-1}$. Furthermore, we have that $s(Z)=s(-Z)$. This contradicts a result of Szab\'o
\cite{refSzabo} and shows that $r=0$. Hence, $\SS(\cdot)$ vanishes identically on
$S^+(T_PM)$. Consequently, $\nabla R=0$ on $T_PM$, see for example
\cite{refGilStav}. \end{proof}

\section{Jordan Osserman manifolds}\label{Sect4}

Let $L$ be the second fundamental form of the hypersurface $(M,g_f)$ and let $R$ be the curvature tensor. We use Lemma
\ref{Lem2.1} and equation (\ref{Eqn1.3}) to see:
\begin{eqnarray*}
&&R(\partial_i^x,\partial_j^x,\partial_k^x,\partial_l^x)=-\textstyle\frac12\{
\partial_i^x\partial_l^x(\partial_j^xf\cdot\partial_k^xf)+
\partial_j^x\partial_k^x(\partial_i^xf\cdot\partial_l^xf)\\&&\qquad\qquad-
\partial_i^x\partial_k^x(\partial_j^xf\cdot\partial_l^xf)-
\partial_j^x\partial_l^x(\partial_i^xf\cdot\partial_k^xf)\}\\
&&\qquad\quad=\partial_i^x\partial_l^xf\cdot\partial_j^x\partial_k^xf-\partial_i^x\partial_k^xf\cdot
\partial_j^x\partial_l^x\\
&&\qquad\quad=L(\partial_i^x,\partial_l^x)L(\partial_j^x,\partial_k^x)-L(\partial_i^x,\partial_k^x)L(\partial_j^x,\partial_l^x).
\end{eqnarray*}
This agrees with the well known formula for the curvature of a hypersurface \cite{refGil}:
\begin{equation}
R(Z_1,Z_2,Z_3,Z_4)=
   L(Z_1,Z_4)L(Z_2,Z_3)-L(Z_1,Z_3)L(Z_2,Z_4).\label{Eqn4.1}
\end{equation}

Assertions (1a) and (2a) of Theorem \ref{Thm1.3} will follow from the following Lemma.

\begin{lemma}\label{Lem4.1} Assume that $L_{\mathcal{X}}$ is non-degenerate on a non-empty connected open subset $\mathcal{O}$ of
$M$.
\begin{enumerate}
\item If $p=2$, then $(\mathcal{O},g_f)$ is Jordan Osserman.
\item If $p\ge3$ and if $L_{\mathcal{X}}$ is definite on $\mathcal{O}$, then $(\mathcal{O},g_f)$ is Jordan Osserman.
\item If $p\ge3$ and if $L_{\mathcal{X}}$ is indefinite on $\mathcal{O}$, then $(\mathcal{O},g_f)$ is neither spacelike Jordan
Osserman nor timelike Jordan Osserman.
\end{enumerate}\end{lemma}

\begin{proof} We use an argument motivated by results of Stavrov \cite{refSt}. Let
$P\in M$ and suppose
$L_{\mathcal X}$ is non-degenerate on
$T_PM$. Let
$Z\in S^\pm(T_PM)$. We decompose
$Z=X+Y$ for
$X\in\mathcal{X}(P)$ and
$Y\in\mathcal{Y}(P)$. Since $(Z,Z)\ne0$ and since $\mathcal{Y}$ is totally isotropic, $X\ne0$. By Lemma
\ref{Lem2.1},
$J(Z)=J(X)$. As $J(X)^2=0$, $\rank(J(X))$ determines  the Jordan normal form of $J(X)$.
Let $0\ne X\in\mathcal{X}(P)$. By equation (\ref{Eqn4.1}), we have:
\begin{eqnarray}
&&(J(X_1)X_2,X_3)
=L(X_1,X_1)L(X_2,X_3)-L(X_1,X_2)L(X_1,X_3).
\label{Eqn4.2}\end{eqnarray}
We have $J(X)\partial_i^y=0$ and
$J(X)X=0$. Thus 
$$\rank(J(X))\le p-1.$$

Suppose first that $L(X,X)\ne0$. We can then choose a basis
$\{X_1,...,X_p\}$ for
$\mathcal{X}(P)$ so $X_1=X$ and so
$L(X_i,X_j)=\varepsilon_i\delta_{ij}$, where $\varepsilon_i\ne0$ for $1\le i\le p$. We use equation (\ref{Eqn4.2}) to show
 $\rank(J(X))=p-1$ by computing:
$$(J(X)X_i,X_j)=\varepsilon_i^2\delta_{ij}\text{ for }i,j\ge2.$$
Suppose next that $L(X,X)=0$. We can then choose a basis so $L(X_1,X_2)=1$ and so $L(X_1,X_i)=0$ for $i\ne2$. We show that
$\rank(J(X))=1$ by computing
$$(J(X)X_1,X_i)=-\delta_{1,i}\quad\text{ and }\quad(J(X)X_i,X_j)=0\text{ for }i\ne 1.$$
Consequently, if $0\ne
X\in\mathcal{X}$, then:
\begin{equation}
\rank(J(X))=\left\{\begin{array}{ll}
p-1&\text{if }\quad(X,X)\ne0,\\
1&\text{if }\quad(X,X)=1.\end{array}\right.\label{Eqn4.3x}\end{equation}

Suppose $L_{\mathcal{X}}$ is definite. Let $X\ne0$. Then $L(X,X)\ne0$, so $\rank(J(X))=p-1$ by equation (\ref{Eqn4.3x}). 
This shows that $(\mathcal{O},g_f)$
is timelike and spacelike Jordan Osserman. If $p=2$, then $p-1=1$. Equation (\ref{Eqn4.3x}) implies $\rank(J(X))=1$
and again
$(\mathcal{O},g_f)$ is timelike and spacelike Jordan Osserman. Finally, if $L_{\mathcal{X}}$ is indefinite and if $p>2$, then
$\rank(J(X))=1$ if $L(X,X)=0$ and
$\rank(J(X))=p-1\ne1$ if $L(X,X)\ne0$. Consequently, $(\mathcal{O},g_f)$ is neither spacelike Jordan Osserman nor timelike Jordan
Osserman. \end{proof}

\section{Jordan IP manifolds}\label{Sect5}
We complete the proof of Theorem \ref{Thm1.3} by proving:

\begin{lemma} Let $p\ge2$.
Assume that $L_{\mathcal{X}}$ is non-degenerate on a non-empty connected open subset $\mathcal{O}$ of $M$.
Then $(\mathcal{O},g_f)$ is:\begin{enumerate}
\item  spacelike Jordan IP and timelike Jordan IP;
\item not mixed Jordan IP.
\end{enumerate}\end{lemma}

\begin{proof} We adopt arguments of \cite{refGiZa} (see Section 5). Let $\{Z_1,Z_2\}$ be an orthonormal
basis for a non-degenerate $2$ plane in $T_PM$ for $P\in\mathcal{O}$. We expand $Z_\nu=X_\nu+Y_\nu$ and use
Lemma \ref{Lem2.1} to see $\RR(\pi)=R(X_1,X_2)$. As $\RR(\pi)^2=0$, $\rank(\RR(\pi))$ determines the Jordan normal form. Equation
(\ref{Eqn4.1}) implies:
\begin{equation}(\RR(\pi)X_3,X_4)=L(X_1,X_4)L(X_2,X_3)-L(X_1,X_3)L(X_2,X_4).\label{Eqn4.4x}\end{equation}

If $\pi$ is spacelike or timelike, then $\pi$ contains no null vectors
and thus $\{X_1,X_2\}$ are linearly independent vectors. We extend this set to a basis $\{X_1,...,X_p\}$ for
$\mathcal{X}(P)$. Since $L_{\mathcal{X}}(P)$ is non-degenerate, we can choose a basis $\{\tilde X_1,...,\tilde X_p\}$
for $\mathcal{X}(P)$ which is dual (with respect to $L$) to the original basis, i.e. $L(X_i,\tilde X_j)=\delta_{ij}$. By
equation (\ref{Eqn4.4x}),
$$(R(X_1,X_2)\tilde X_i,\tilde X_j)=\left\{\begin{array}{rl}
1&\text{ if }i=2,j=1,\\
-1&\text{ if }i=1,j=2,\\
0&\text{ otherwise}.\end{array}\right.$$
It now follows that $\rank(\RR(\pi))=\rank(R(X_1,X_2))=2$. Thus $(\mathcal{O},g_f)$ is spacelike Jordan IP and timelike Jordan IP.

To see that $(\mathcal{O},g_f)$ is not mixed IP, we consider the following $2$ planes:
$$\pi_1:=\Pspan\{\partial_1^y,\partial_1^x\}\quad\text{and}\quad
\pi_2(\varepsilon):=\Pspan\{\varepsilon^{-1}\partial_1^y+\varepsilon\partial_1^x,
-\varepsilon^{-1}\partial_2^y+\varepsilon\partial_2^x\},$$
respectively, where $\varepsilon$ is a real parameter.
The matrices giving the induced inner products on $\pi_1$ and $\pi_2(\varepsilon)$ are given by:
$$A_1:=\left(\begin{array}{ll}0&1\\1&\varrho_{11}\end{array}\right)\quad\text{and}\quad
A_2:=\left(\begin{array}{ll}2+\varepsilon^2\varrho_{11}&\varepsilon^2\varrho_{12}\\
\varepsilon^2\varrho_{12}&-2-2\varepsilon^2\varrho_{22}
\end{array}\right),$$
where $\varrho_{ij}:=(\partial_i^x,\partial_j^x)$.
Since $\det(A_1)=-1$ and $\det(A_2)=-4+O(\varepsilon^2)$, $\pi_1$ and $\pi_2(\varepsilon)$ are mixed $2$ planes for $\varepsilon$
small. Since
$\RR(\pi_1)=0$ and $\RR(\pi_2(\varepsilon))=c(\varepsilon)R(\partial_1^x,\partial_2^x)\ne0$, $\RR(\pi_1)$ and
$\RR(\pi_2(\varepsilon))$ are not Jordan equivalent and hence
$(\mathcal{O},g_f)$ is not mixed Jordan IP. \end{proof}

\section{Manifolds of signature $p\ne q$}\label{Sect6}

\begin{proof}[Proof of Theorem \ref{Thm1.5}]
Let $N$ be the isometric product of $\mathbb{R}^{(a,b)}$ with $(M,g_f)$ this has signature $(p+a,q+b)$. Let
$R^N$ and
$R^M$ be the curvature tensors on $N$ and $M$ respectively. Let
$U_\nu$ be tangent vectors on $N$. We decompose $U_\nu=W_\nu+Z_\nu$, where $W_\nu$ is tangent to
$\mathbb{R}^{(u,v)}$ and
$Z_\nu$ is tangent to $M$. Since
\begin{eqnarray*}
&&R^N(U_1,U_2)U_3=R^M(Z_1,Z_2)Z_3,\text{ and }\\
&&\nabla_{U_1}^NR^N(U_2,U_3)=\nabla_{Z_1}^MR^M(Z_2,BZ_3),\end{eqnarray*}
$J^N(U)=J^N(Z)$, $\SS^N(U)=\SS^M(Z)$, and $R^N(U_1,U_2)=R^M(Z_1,Z_2)$. We use
(\ref{Eqn2.5}) to see that
$(N,g_N)$ is nilpotent Osserman Szab\'o IP; this proves assertion (1); assertions (2) and (3) follow from the corresponding
assertions for $(M,g_f)$.

Suppose that $b=0$. Let $0\ne U\in TN$ be spacelike. Expand $U=W_Z$ and $Z=X+Y$. If
$X=0$, then
$Z\in\mathcal{Y}$ so $(Z,Z)=(U,U)\le 0$, which is false. Thus
$X\ne0$ and by equation (\ref{Eqn4.3x})
$$\rank(J(U))=\rank(J(X))=p-1.$$ 
Thus $(N,g_N)$ is spacelike Jordan Osserman. One shows similarly that $(N,g_N)$ is spacelike Jordan IP. This proves assertion (4b);
assertion (5a) follows similarly.

Suppose $b>0$. We can choose $0\ne W\in
T(\mathbb{R}^{(a,b)})$ spacelike. Then we have that $\rank(J(W))=0$. We can choose $0\ne Z\in TM$ spacelike so $\rank(J(Z))=p-1$.
Thus
$(N,g_N)$ is not timelike Jordan Osserman. Similarly, we may show that $(N,g_N)$ is not timelike IP. This proves assertions (5b)
and (6b); the proof of assertion (6a) is similar.
\end{proof}

\bigbreak\noindent{\bf Acknowledgments:} Research of P. Gilkey partially supported by the NSF (USA) and
the MPI (Leipzig); research of R. Ivanova and T. Zhang partially supported by the NSF (USA).

\end{document}